\author{Rebecca M. Crossley$^{1}$, Carles Falc\'o$^{1}$, and Ruth E. Baker$^{1}$}
\date{}
\title{ {An optimal control approach to nonlinear wave speed selection in reaction--diffusion equations}}

\documentclass[onecolumn]{custom_template}

\usepackage{natbib}

\usepackage{amsmath,amssymb}
\usepackage{graphicx}
\usepackage{xcolor}
\usepackage{amsfonts}

\begin{document}


\maketitle

\begin{center}
{\vspace{-.85cm}$^{1}$\emph{Mathematical Institute, University of Oxford, Oxford, United Kingdom OX2 6GG }}
\end{center}


\begin{abstract}
Travelling wave solutions of reaction–diffusion equations are widely used to model the spatial spread of populations and other phenomena in biology and physics. 
In this article, we reinterpret the classical variational principle approach through an optimal control formulation, in order to obtain a lower bound on the invasion speed of travelling wave solutions in systems of nonlinear partial differential equations. 
We begin by analysing single--species models, where the evolution of the density is governed by a scalar equation with a density-dependent diffusion term and a nonlinear reaction term. 
We show that for any admissible test function, maximising with respect to the parameter of interest yields a bound on the travelling wave speed. 
We apply this framework to several examples, including the porous--Fisher equation, and examine when nonlinear selection mechanisms dominate over the classical linear marginal stability criterion.
Extending this approach, we then consider multi-species systems of reaction--diffusion equations and, reframed as Pontryagin-type optimality systems, we derive analogous bounds on the travelling wave speed using a variational framework under weak coupling. 
Finally, we employ numerical simulations to confirm the accuracy of the predicted wave speeds across a range of illustrative examples.
\end{abstract}

\section{Introduction}
Reaction--diffusion equations have long served as a powerful and widely-used mathematical framework for modelling the collective movement of populations of individuals. 
These models typically take the form of partial differential equations (PDEs) with density-dependent reaction (growth) and diffusion terms, enabling a detailed description of both the spatial and temporal evolution of a population across a  variety of biological and ecological contexts~\cite{murray_mathematical_1993}.

A hallmark feature of such models is the emergence of travelling wave solutions that maintain a fixed profile while propagating at a constant speed. 
These travelling waves have become central to our theoretical understanding of diverse biological processes, from wound healing to tumour invasion, and species migration~\cite{simpson2007cell, trewenack2009traveling,gerlee2016travelling, browning2019bayesian, falco2023quantifying}.
In many of these settings, the wave connects an unstable steady state (\textit{e.g.}, uncolonised space) to a stable steady state (\textit{e.g.}, an established population), 
even in the absence of classical linear instabilities. 
This is unlike systems such as crystal growth, nerve impulses, or fluid flow, where instabilities in the linear system are key drivers of wave propagation. 

One of the seminal models exhibiting travelling wave solutions is the Fisher--KPP model~\cite{fisher_wave_1937, kolmogorov1937study, murray_mathematical_1993}. 
Along with its many extensions, the Fisher--KPP model clearly illustrates how changes in the model parameters shape the speed and form of the invading wave~\cite{crossley2023travelling, fisher_wave_1937, kolmogorov1937study, murray_mathematical_1993}.
Considerable analytical effort has focused on understanding the influence of the various nonlinearities in single--species reaction--diffusion equations, with recent work yielding increasingly refined existence and predicted wave speed results~\cite{bouin2014travelling, colson2021travelling, de1998travelling, el2021travelling, FALCO2024109209,gerlee2016travelling,hadeler1975travelling,  marchant2001travelling,  sanchez1994travelling, stepien2018traveling,yue2025traveling}. 

Several complementary perspectives have emerged to explain how travelling wave speeds are selected, each with its own strengths and limitations. 
A foundational viewpoint uses the linear marginal stability criterion, which identifies the travelling wave with the steepest permissible exponential decay into the unstable state~\cite{van1989front}. 
While this criterion provides an elegant and computable prediction, it relies on linearisation around the leading edge and therefore cannot account for cases where nonlinear interactions---either in the bulk or in the front---influence the selected speed. 
To address this limitation, Benguria and Depassier introduced a variational method for wave speed selection, initially for nonlinear reaction terms~\cite{benguria1996speed, benguria1996variational} and later extended to incorporate nonlinear diffusion~\cite{BENGURIA2022112668}. 
Their approach generates rigorous upper and lower bounds on the admissible wave speeds and, crucially, provides a selection criterion identifying when the minimal admissible speed is indeed the realised one~\cite{benguria1994validity}. 
Hadeler and Rothe~\cite{hadeler1975travelling} clarified that, under suitable conditions, the asymptotic propagation speed can coincide with the linearly predicted value, even when the variational upper and lower bounds do not match. 
These results highlight cases where the linear and nonlinear speed selection criteria actually coincide.
Though broadly undeveloped thereafter, more recent studies, such as that of Stokes et al.~\cite{stokes2024speed}, have since extended these ideas by systematically analysing how density-dependent diffusion shapes both travelling wave existence and speed, thereby revealing parameter regimes where nonlinear diffusion either reinforces or invalidates predictions from earlier linear or variational theories.


Taken together, these developments establish a rich set of tools for analysing travelling waves in single--species systems, yet they also expose a notable gap: most classical speed-selection principles, including both the marginal stability criterion and the Benguria--Depassier variational framework, are formulated for scalar equations. 
However, many biological invasion processes fundamentally involve interactions between multiple populations.
Extending the theory to determine travelling wave speed selection in these multi-species systems is therefore essential for connecting mechanistic models of collective migration with the analytical insights that have proved so valuable in the single--species setting. 

We analyse a range of mathematical models with different combinations of density-dependent reaction and diffusion terms. 
Analytical techniques are used where feasible, supported by numerical simulations implemented using a finite-difference scheme described in~\cite{crossley2023travelling}, to illustrate the existence and selection of travelling waves across different parameter regimes. 
The main contributions of this paper are threefold. 
First, we reinterpret the classical variational principle for travelling wave speed selection as an optimal control problem, clarifying its structure and limitations. 
Second, we extend this formulation to multi-species reaction--diffusion systems, where classical variational principles are no longer directly applicable. 
Third, we derive a tractable variational approximation of the resulting optimality system that yields sharp lower bounds on travelling wave speeds in a range of biologically motivated models, including systems exhibiting nonlinear speed selection.

The remainder of the paper is organised as follows.
First we present the optimal control approach to nonlinear wave speed selection for general multispecies models (Sec.~\ref{sec:OC}). 
Next, in Sec.~\ref{1-vp} we revisit the classical variational principle for a single--species and provide a sufficient criterion for nonlinear wave speed selection in terms of the functional forms of the diffusion and reaction terms.
This scalar setting provides a benchmark against which we can assess how speed-selection mechanisms change when additional interactions are incorporated.
We then apply the optimal control framework to a number of examples in the literature including for models involving two interacting species (Sec.~\ref{sec:MS}). 
We show that this approach successfully predicts a minimum travelling wave speed in a number of examples, including cell invasion and differentiation. 
We conclude in Sec.~\ref{sec:conc} with a brief discussion of these results and suggest some potential avenues for future research.

\section{An optimal control approach for multispecies models}\label{sec:OC}

We consider reaction--diffusion systems for $n\geq1$ species
$\boldsymbol{\rho}=(\rho_1,\Tilde{\boldsymbol{\rho}})^\top$, with
$\Tilde{\boldsymbol{\rho}}=(\rho_2,\ldots,\rho_n)^\top$, of the form
\begin{align}
    \partial_t\rho_1 & = \partial_x\!\left(D(\boldsymbol{\rho})\,\partial_x\rho_1\right) + f(\boldsymbol{\rho}), \\
    \partial_t\Tilde{\boldsymbol{\rho}} & = \mathbf{g}(\boldsymbol{\rho}),
\end{align}
where $\rho_i = \rho_i(x,t)$ with $x\in\mathbb{R},\,t\geq 0$, and $D,f:\mathbb{R}^n\to\mathbb{R}$, $\mathbf{g}:\mathbb{R}^n\to\mathbb{R}^{n-1}$. 
This type of model considers a single motile population, $\rho_1$, with the rest of species evolving via reaction terms---these are often used, for instance, in systems describing cell invasion into the extracellular matrix~\cite{browning2019bayesian, colson2021travelling, crossley2023travelling}. 
We focus on systems admitting constant speed ($c\geq 0$), constant profile travelling wave solutions $\boldsymbol{\rho}(x,t)=\boldsymbol{u}(z)$,
$z=x-ct$, where $\boldsymbol{u}=(u_1,\Tilde{\boldsymbol{u}})^\top$ and
$\Tilde{\boldsymbol{u}}=(u_2,\ldots,u_n)^\top$.
We assume that $u_1$ is
non-increasing and, after nondimensionalisation,
\begin{equation*}
    \begin{cases}
        u_1(z)\to1 & \text{as } z\to-\infty, \\
        u_1(z)\to0 & \text{as } z\to+\infty .
    \end{cases}
\end{equation*}
The travelling wave solution also fixes boundary conditions for $\Tilde{\boldsymbol{u}}$ for $z\rightarrow\pm\infty$.

In travelling wave coordinates, and defining
$v=-\mathrm{d}u_1/\mathrm{d}z\geq0$, the system becomes
\begin{align*}
    \frac{\mathrm{d}}{\mathrm{d}z}\!\left(D(\boldsymbol{u})\,v\right)
    & = -cv + f(\boldsymbol{u}), \\
    c\,\frac{\mathrm{d}\Tilde{\boldsymbol{u}}}{\mathrm{d}z}
    & = -\mathbf{g}(\boldsymbol{u}) .
\end{align*}
Using the monotonicity of $u_1$, we can rewrite the system in terms of $u_1$,
with $\mathrm{d}/\mathrm{d}z=-v(u_1)\,\mathrm{d}/\mathrm{d}u_1$, yielding
\begin{align}
v(u_1)\,\frac{\mathrm{d}}{\mathrm{d}u_1}\!\left(D(\boldsymbol{u})\,v(u_1)\right)
    & = cv(u_1) - f(\boldsymbol{u}), \label{eq:v thing}\\
cv(u_1)\frac{\mathrm{d}\Tilde{\boldsymbol{u}}}{\mathrm{d}u_1}
    & = \mathbf{g}(\boldsymbol{u}) .
    \label{eq: odes tw}
\end{align}

To obtain a variational characterisation of the wave speed $c$, we
introduce a non-negative, non-increasing test
function $\varphi(u_1)$, with $\phi(u_1)=-\varphi'(u_1)\geq0$ and $\varphi(1)=0$.
We note that for systems ($n>1$), the remaining equations (Eq.~\eqref{eq: odes tw}) act as
dynamical constraints. 
Multiplying Eq.~\eqref{eq:v thing} by $D(\boldsymbol{u})\varphi(u_1)$ and applying integration by parts gives
\begin{align*}
\mathcal{I}_c[v\,;\Tilde{\boldsymbol{u}}]
:&= \int_0^1 \!\left[
    -\tfrac12 \phi(u_1)\,D(\boldsymbol{u})^2\,v(u_1)^2
    + cv(u_1)\,D(\boldsymbol{u})\,\varphi(u_1)
    - D(\boldsymbol{u})\,f(\boldsymbol{u})\,\varphi(u_1)
\right]\mathrm{d}u_1 \\
&= 0 .
\end{align*}

The quadratic dependence of $\mathcal{I}_c$ on $v$ suggests studying the associated variational problem, with the aim of deriving bounds that may provide direct estimates of the wave speed, $c$, via
\begin{equation*}
    0 = \mathcal{I}_c[v\,;\Tilde{\boldsymbol{u}}] \leq \mathcal{I}_c[v^*\,;\Tilde{\boldsymbol{u}^*}]\,,
\end{equation*} for an optimal pair $(v^*,\Tilde{\boldsymbol{u}}^*)$. 
In systems, however, $v$
appears implicitly through Eq.~\eqref{eq: odes tw}. We therefore consider the
optimisation problem
\begin{equation}
    v^*= \arg\max_v \; \mathcal{I}_c\!\left[v\,;\Tilde{\boldsymbol{u}}\right],\qquad \text{subject to}\;\; cv(u_1)\,\frac{\mathrm{d}\Tilde{\boldsymbol{u}}}{\mathrm{d}u_1}
    = \mathbf{g}(\boldsymbol{u}) .
\end{equation}

This is a standard optimal control problem, which we solve using
Pontryagin's maximum principle~\cite{lenhart2007optimal}. 
Introducing the Hamiltonian defined by
\begin{align*}
    \mathcal{H}[u_1,v\,;\Tilde{\boldsymbol{u}},\boldsymbol{\omega}]
    &=
    -\tfrac12 \phi(u_1)\,D(\boldsymbol{u})^2\,v(u_1)^2
    + c\,v(u_1)\,D(\boldsymbol{u})\,\varphi(u_1)
    \\ & \qquad\qquad\qquad\qquad- D(\boldsymbol{u})\,f(\boldsymbol{u})\,\varphi(u_1)
    + \frac{\boldsymbol{\omega}^\top\cdot\mathbf{g}(\boldsymbol{u})}{cv(u_1)} ,
\end{align*}
where $\boldsymbol{\omega}=(\omega_2,\ldots,\omega_n)^\top$ are co-state
variables. Pontryagin's principle implies that maximisers $v^*$ are
pointwise maximisers of $\mathcal{H}$. When the Hamiltonian is concave down in
$v$, these satisfy
\begin{equation*}
    \frac{\partial \mathcal{H}}{\partial v} = -\phi(u_1)\,D(\boldsymbol{u})\,v(u_1)+c\,D(\boldsymbol{u})\,\varphi(u_1) - \frac{\boldsymbol{\omega}^\top\cdot\mathbf{g}(\boldsymbol{u})}{cv(u_1)^2} =  0\,,
\end{equation*}
or equivalently
\begin{equation}
    -\phi(u_1)\,D(\boldsymbol{u})^2\,v(u_1)^3+c\, D(\boldsymbol{u})\,\varphi(u_1)\,v(u_1)^2 - \frac{\boldsymbol{\omega}^\top\cdot\mathbf{g}(\boldsymbol{u})}{c} =  0\,. \label{eq:v3}
\end{equation}
Eq.~\eqref{eq:v3} implicitly determines the optimiser $v^{*}$ as the positive real root of the cubic equation.
Furthermore, the co-state variables evolve according to
\begin{align*}
    \frac{\mathrm{d}\boldsymbol{\omega}}{\mathrm{d}u_1}
   &  = -\nabla_{\Tilde{\boldsymbol{u}}}\mathcal{H} \\
   & = \left[\phi(u_1)\,D(\boldsymbol{u})\,v(u_1)^2- cv(u_1)\,\varphi(u_1)\right]\nabla_{\Tilde{\boldsymbol{u}}} D(\boldsymbol{u}) \\ &\qquad\qquad\qquad\qquad+\nabla_{\Tilde{\boldsymbol{u}}}(D(\boldsymbol{u})f(\boldsymbol{u}))\,\varphi(u_1)-\frac{\boldsymbol{\omega}^\top \nabla_{\Tilde{\boldsymbol{u}}}\mathbf{g}(\boldsymbol{u})}{cv(u_1)}.
\end{align*}
No transversality conditions are imposed on $\boldsymbol{\omega}$, since
travelling wave solutions prescribe boundary conditions for the state variables 
$\Tilde{\boldsymbol{u}}$ at $u_1=0$ and $u_1=1$.

In this next section, we analyse this problem for the specific case of one-species models. 
In doing so, we unveil a relationship to the variational principle for one-species models, previously analysed in~\cite{benguria1994validity, benguria1996speed, benguria1996variational, stokes2024speed}, before applying the optimal control approach to a number of new examples that illustrate why the optimisation yields useful bounds on the travelling wave speed $c$. 

\section{One-species models: recovering the variational principle}\label{1-vp}
In cases where there is only one species, $\boldsymbol{\rho} = (\rho_1)$, the dynamical constraint imposed by Eq.~\eqref{eq: odes tw} vanishes. 
Here, we require that $f(\rho_1)$ has two equilibrium points, we choose these without loss of generality at $\rho_1 = 0$, and $\rho_1 = 1$, \textit{i.e.}, $f(0) = f(1) = 0$.
In this case, maximising the functional $\mathcal{I}_c$ relies on standard variational arguments. 
In particular, omitting indices for simplicity, we have
\begin{equation}
    \mathcal{I}_c[v]
:= \int_0^1 \!\left[
    -\tfrac12 \phi(u)\,D(u)^2\,v(u)^2
    + c\,v(u)\,D(u)\,\varphi(u)
    - D(u)\,f(u)\,\varphi(u)
\right]\mathrm{d}u
= 0 .\label{eq:func_1}
\end{equation}
 Interpreting the integrand in Eq.~\eqref{eq:func_1} as a quadratic polynomial in $v$ reveals a natural inequality:
\begin{equation}
-\frac{1}{2}\phi(u)\,D(u)^2\,v(u)^2+c\,v(u)\,D(u)\,\varphi(u) \leq \frac{c^2\varphi(u)^2}{2\phi(u)}, \label{eq:bound_1}
\end{equation}
that transforms the integral identity into a bound that holds for all admissible test functions. 
Equality in Eq.~\eqref{eq:bound_1} is only obtained when  
\begin{equation}
    v(u) = v^*(u) := \frac{c \varphi(u)}{D(u)\phi(u)}\,. \label{eq:vstar_VP_1}
\end{equation}
Inserting the bound in Eq.~\eqref{eq:bound_1} into Eq.~\eqref{eq:func_1} yields a family of lower bounds on the wave speed, indexed by the choice of smooth test function $\varphi(u)$ with $\varphi'(u) = - \phi(u)$:
\begin{equation}
    \frac{1}{2}c^2 \geq \frac{\int_0^1 D(u)f(u)\varphi(u)\,\mathrm{d}u}{\int_0^1 \varphi(u)^2/\phi(u)\,\mathrm{d}u},
    \label{eq: speed bound one species}
\end{equation}

To make the bound explicit, we now select a convenient one-parameter family of test functions, $\varphi_\beta (u) = ((1-u)/u)^\beta$, for $\beta \in[0,2)$. 
This choice, previously used in~\cite{benguria1994validity, stokes2024speed}; satisfies the condition that the integral remains finite and boundary terms vanish in the limiting sense. 
As such, it balances analytical tractability with the ability to capture the asymptotic behaviour near the steady states.
Using this test function, we find
\begin{equation}
    \frac{1}{2}c^2\geq \sup_{\beta\in [0,2)} F(\beta),\quad\mbox{with } F(\beta) = \beta\,\frac{\int_0^1 D(u)f(u) u^{-\beta} (1-u)^{\beta} \,\mathrm{d}u}{{B}(2-\beta, 2+\beta)}\,,
    \label{eq: bound one species F}
\end{equation}
where $\mathrm{B}$ is the Beta function.

To confirm that the inequality in Eq.~\eqref{eq: bound one species F} defines a genuine variational principle, we must also show that there exists a test function, $\varphi=\hat{\varphi}$ (and correspondingly, $\phi=\hat{\phi}$) for which the bound is attained. 
This requires analysing the asymptotic behaviour of the travelling wave profile near $u=0$ and $u=1$ and verifying that the corresponding integrals converge.
Equality in Eq.~\eqref{eq: speed bound one species} is attained (as per Eq.~\eqref{eq:vstar_VP_1}) when 
\begin{equation}
    v(u)=\dfrac{c\hat{\varphi}(u)}{D(u)\hat{\phi}(u)} = -\dfrac{c\hat{\varphi}(u)}{D(u)\hat{\varphi}'(u)}, \label{eq: equal 1}
\end{equation}
where $'=\mathrm{d}/\mathrm{d}u.$
Eq.~\eqref{eq: equal 1} can be integrated to obtain
\begin{equation}
    \hat{\varphi}(u) = \exp\Bigg\{-c\int_{\bar{u}}^{u}\dfrac{\mathrm{d}q}{D(q)v(q)}\Bigg\}, \label{eq:hat_varphi}
\end{equation}
for some $\bar{u}\in(0,1).$
Clearly $\hat{\varphi}>0$ and is monotonically decreasing.
We now make the change of variables using $p(u)=D(u)v(u).$
Eq.~\eqref{eq:v thing} becomes 
\begin{equation}
    p(u)\dfrac{\mathrm{d}}{\mathrm{d}u}p(u)-cp(u)+D(u)f(u)=0.\label{eq:p}
\end{equation}
In order to ensure convergence of $\hat{\varphi}(u)$ across all of $u\in[0,1]$, we must then analyse the asymptotic behaviour of Eq.~\eqref{eq:p} as $u\to0^{+}$ and $u\to1^{-}$~\cite{benguria2004minimal}.

Begin by considering the linear analysis around the endpoint $u=0.$ 
Near $u=0$, we can write $f(u)\sim f'(0)u$ and $D(u)\sim D(0)$ (if $D(0)=0$ or $f'(0) = 0$, then the marginal stability analysis predicts a zero wave speed), so that by employing the ansatz $p(u)=\alpha u$, where $\alpha$ is the larger root of $m^2-cm+D(0)f'(0)=0$, we find that 
\begin{equation*}
    c\geq 2\sqrt{D(0)f'(0)}=c_L,
\end{equation*}
in order for $\alpha$ to be real, and hence avoid oscillatory solutions that predict negative densities. 
Hereon in, $c_L$ is the travelling wave speed predicted by linear stability analysis.
Solving Eq.~\eqref{eq:hat_varphi} for $p(u)=\alpha u$ and $\bar{u}=u(z\to-\infty)=1$ gives $\hat{\varphi}=u^{-c/\alpha}$, which ensures that all integrals in Eq.~\eqref{eq: speed bound one species} converge as long as $c/\alpha<2.$

Now consider the linear analysis around the endpoint $u=1$. 
Define the new variable, $s$, as $s=1-u$. 
Using Taylor series expansions, $f(u)\sim f'(1)(u-1)=-|f'(1)|s$ and $D(u)\sim D(1).$
Setting $p(u)=r(1-u)=rs$ and employing this ansatz in Eq.~\eqref{eq:p} then, following the same methods as before, we find that for $\bar{u}=1$, $\hat{\varphi}=(1-u)^{c/r}$ with $r$ as the positive root of the auxiliary equation $r^2+cr+f'(1)D(1)=0$, giving
\begin{equation*}
    r = \dfrac{-c\pm\sqrt{c^2-4D(1)f'(1)}}{2}.
\end{equation*}
As such, the integrals in the functional~\eqref{eq: bound one species F} must converge at $u=1$ because $c/r>0$ and does not restrict the speed. 
This is sufficient to show convergence of all of the integrals across the entire domain $[0,1]$ for an attainable test function $\hat{\varphi}$, so that Eq.~\eqref{eq: bound one species F} does in fact define a variational principle.

\subsection{A criterion for nonlinear wave speed selection}\label{sec:criterion}
The classical variational framework provides lower bounds on the admissible travelling wave speed but does not, in general, indicate whether the minimal speed is selected by linear leading-edge dynamics or by nonlinear effects in the interior of the wave.
In this section, we extract from the variational formulation a simple diagnostic criterion that predicts when nonlinear (``pushed'') wave speed selection occurs, \textit{i.e.}, when $c>c_L$. 
The criterion is based on the behaviour of the variational functional $F(\beta)$ (Eq.~\eqref{eq: speed bound one species}) near the boundary $\beta=2$, which corresponds to the pulled wave regime where the travelling wave speed is governed by the leading edge dynamics. 
The analysis that follows builds on earlier variational studies~\cite{benguria1996speed, stokes2024speed}, but is reformulated here to yield an explicit condition. 

If $D$ is continuous near the origin, $f\in C^1$ near the origin, and $D(u)f(u)(1-u)^\beta$ is bounded, then the dominant contribution of the integral in the numerator of $F(\beta)$ as $\beta\nearrow2$ arises from a neighbourhood of $u=0$.
In particular, if $D(0)f'(0)\neq0$, then
\begin{equation*}
    \int_0^1 D(u)f(u) u^{-\beta} (1-u)^{\beta} \,\mathrm{d}u = D(0)f'(0)\int_0^1 u^{-\beta + 1}\,\mathrm{d}u + \mathcal{O}(1) = \frac{D(0)f'(0)}{2- \beta} + \mathcal{O}(1)\,.
\end{equation*}
Similarly, the Beta function satisfies
$\mathrm{B}(2-\beta,2+\beta)=(2-\beta)^{-1}-11/6+\mathcal{O}(2-\beta)$, and hence
\begin{equation*}
    \lim_{\beta\nearrow 2} F(\beta) = 2D(0)f'(0)\,,
\end{equation*}
which recovers the linear marginal stability bound
$c\ge c_L:=2\sqrt{f'(0)D(0)}$.

Since the minimal admissible wave speed is determined by
\[
c^2 = 2\sup_{\beta\in[0,2)} F(\beta),
\]
we see that nonlinear (pushed) wave speed selection occurs precisely when the supremum of $F(\beta)$ is attained at an interior point $\beta<2$. 
A sufficient local condition for this is that $F(\beta)$ increases as $\beta$ decreases from $\beta=2$, namely
\begin{equation}
\left.\frac{\mathrm{d}}{\mathrm{d}\beta}F(\beta)\right|_{\beta=2^-}<0.\label{eq:beta_criteria}
\end{equation}
We therefore study the sign of $F'(\beta)$ as $\beta\nearrow2$, assuming that $F''(\beta)\le0$ on $[0,2)$ so that $F$ is concave.

This criterion (Eq.~\eqref{eq:beta_criteria}) immediately explains several known regimes.
If $f(u)\ge0$ for all $u\in[0,1]$ and $f'(0)D(0)=0$, then $F(\beta)\ge0$ and the wave speed is necessarily selected nonlinearly. 
This includes Porous--Fisher-type models with $D(u)\sim u^m$ as $u\to0$, as well as models with Allee-type reaction terms, where $f(u)= u(u-a)(1-a)$ is negative at low densities.
When $D(0)f'(0)\neq0$, the sign of $F'(2^-)$ determines whether nonlinearities in the diffusion or reaction terms shift the selected speed away from $c_L$.
If $F'(\beta)<0$, then the speed is nonlinearly selected, as we can find a lower bound for the speed on the interval $\beta\in(0,2)$, which exceeds $c_L$. 
If we set $R(u)$ through the relationship $f(u) = uR(u)$, we can define
\begin{equation*}
    N(\beta) = \int_0^1 D(u)f(u) u^{-\beta} (1-u)^{\beta} \,\mathrm{d}u = \int_0^1 D(u)R(u) u^{1-\beta} (1-u)^{\beta} \,\mathrm{d}u\,, 
\end{equation*}
with $R(0)\neq 0$. We further assume that 
\begin{equation*}
    \int_0^1\frac{D(u)R(u)-R(0)D(0)}{u}\,\mathrm{d}u < \infty\,.
\end{equation*}
In particular, this means that $D(u)R(u) - D(0)R(0) =\mathcal{O}(u^a)$ for any $a > 0$. In this case
\begin{align*}
    N(\beta) & = D(0)R(0)\mathrm{B}(2-\beta, 1+\beta) + \int_0^1\frac{D(u)R(u)-R(0)D(0)}{u^{\beta-1}}(1-u)^\beta\,\mathrm{d}u
    \\
    & = \frac{D(0)R(0)}{2-\beta} - \frac{3D(0)R(0)}{2} + \int_0^1\frac{D(u)R(u)-R(0)D(0)}{u^{\beta-1}}(1-u)^\beta\,\mathrm{d}u + \mathcal{O}(2-\beta)\,.
\end{align*}
Using this expression for $N(\beta)$ together with the asymptotic expansion of the Beta function, we obtain
\begin{align*}
    F(\beta) &= \underbrace{2D(0)R(0)}_{F(2)} \nonumber \\ &\qquad - \underbrace{\frac{2}{3}\left[ \frac{1}{2}D(0)R(0) - 3 \int_0^1\frac{D(u)R(u)-R(0)D(0)}{u}(1-u)^2\,\mathrm{d}u \right]}_{F'(2)}(2-\beta)\\&\qquad\qquad+ \mathcal{O}((2-\beta)^2)\,.
\end{align*}
Hence we obtain a sufficient condition for  nonlinear speed selection (\textit{i.e.} $F'(2) < 0$)
\begin{equation}
   \int_0^1\frac{D(u)R(u)-R(0)D(0)}{u}(1-u)^2\,\mathrm{d}u > \frac{1}{6}D(0)R(0)\,,\label{eq: criterion one species}
\end{equation}
where we note that $R(0) = f'(0)$. For instance, if $D(u) = u + \delta$ and $f(u) = u(1-u)$, we have $R(u) = 1-u$ and $D(u)R(u) - D(0)R(0) = u(1-u-\delta)$. 
Eq.~\eqref{eq: criterion one species} predicts a nonlinear contribution to the selected wave speed when $\delta < {1}/{2}$, agreeing with previous results~\cite{stokes2024speed}.

The criterion given in Eq.~\eqref{eq: criterion one species} reflects a competition between leading edge linearisation and nonlinear contributions in the wave interior, and this mechanism persists in multi-species systems under weak or smooth coupling---see Sec.~\ref{sec:MS}.
The remainder of this section focuses on applying the variational approach to a number of one-species models adapted from the literature to improve the estimated lower bound on the travelling wave speed.

\subsection{A model of epidermal wound healing}
\begin{figure}[htbp]
    \centering
    \includegraphics[width=0.75\linewidth]{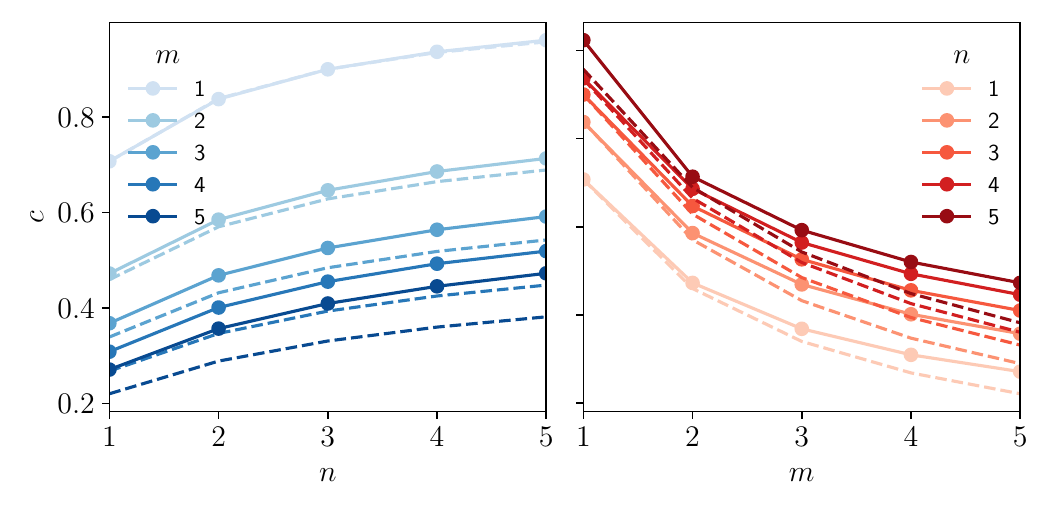}
    \caption{Plot of the numerically estimated minimum travelling wave speed (solid lines) and the minimum travelling wave speed of Eq.~\eqref{eq:wound} predicted using the variational principle with $F(\beta)$ as given in Eq.~\eqref{eq:F_wound} (dashed lines) as a function of $n$ (left-hand side) and $m$ (right-hand side). The numerically estimated travelling wave speed is obtained by tracing the point $X(t)$ such that $u(X(t), t)=0.1$.}
    \label{fig:wound-healing}
\end{figure}

First, we consider a generalised, non-dimensional equation which can be used to study wound healing and other applications~\cite{falco2023quantifying,liu2023parameter,sherratt1990models}. The governing reaction--diffusion equation takes the form
\begin{equation}
    \partial_t\rho=\partial_x(\rho^m\partial_x\rho)+\rho(1-\rho^n), \label{eq:wound}
\end{equation}
for positive integers $m,\,n$.
This is an example of the well-studied Porous--Fisher equation, with a Richards' growth term~\cite{de1998travelling, witelski1995merging}.
We apply the variational principle to estimate the wave speed, using $D(u) = u^{m}$, and $f(u) = u(1-u^n)$. 
In this case, $F(\beta)$ can be  calculated explicitly as
\begin{align}
    F(\beta) &= \dfrac{\beta}{\mathrm{B}(2-\beta, 2+ \beta)} \Bigg[\mathrm{B}(m-\beta+2, \beta+1)-\mathrm{B}(m+n-\beta+1, \beta+1) \Bigg] \nonumber \\
    &= \dfrac{6\beta}{(1+\beta)\Gamma(2-\beta)} \Bigg[\dfrac{\Gamma(m-\beta+2)}{\Gamma(m+3)}-\dfrac{\Gamma(m+n-\beta+2)}{\Gamma(m+n+3)}\Bigg]. \label{eq:F_wound} 
\end{align}
Taking the expression for $F(\beta)$ from Eq.~\eqref{eq:F_wound} and evaluating the supremum in Eq.~\eqref{eq: bound one species F} gives a new lower bound on the travelling wave speed, $c$, which cannot be evaluated explicitly, but is plotted numerically in Fig.~\ref{fig:wound-healing}.

Several key features are apparent in Fig.~\ref{fig:wound-healing}.
First, across all parameter regimes shown, the variational principle provides a strict lower bound, with the dashed curves always lying below the numerically observed minimum travelling wave speeds. 
The bound is tightest when the diffusivity is closest to linear, \textit{i.e.}, for small values of $m$, and when the reaction term is closest to logistic growth, \textit{i.e.}, for small values of $n$. 
In these regimes, the wave profile behaves more similarly to the Fisher--KPP equation, for which the variational bound is known to be exact.
Furthermore, as $m$ increases, the diffusion exhibits a stronger degeneracy, leading to a sharper wavefront.
In this limit, the variational bound underestimates the travelling wave speed more strongly as the impact of the diffusion term increases. 
A similar trend occurs when $n$ increases: the reaction term becomes more strongly nonlinear, moving the system away from the Fisher--KPP regime in which the variational bound is sharp.
In both cases, the gap between the dashed and solid curves reflects precisely the departure from the classical marginal-stability scenario that underpins the variational method.

These trends motivate a closer examination of the special case $n=1$ (\textit{i.e.}, the Porous--Fisher equation), for which the variational expression simplifies substantially and exact analytical results are available for several integer values of $m$.
When $n=1$, we use Eq.~\eqref{eq: bound one species F} to find 
\begin{equation}
    F(\beta) = \beta\frac{\mathrm{B}(m+2-\beta,2+\beta)}{\mathrm{B}(2-\beta, 2+ \beta)} = \beta \frac{\Gamma(4)}{\Gamma(4+m)}\frac{\Gamma(2-\beta + m)}{\Gamma(2-\beta)}, \label{eq:F_n1}
\end{equation}
which can be simplified for small integer values of $m$, as $F(\beta)$ is a low-order polynomial in $\beta$. 
For instance, we recover the known values of the travelling wave speed when $m = 0$ and $m = 1$. 
For $m=0$, Eq.~\eqref{eq:F_n1} reduces to $F(\beta)=\beta$, and thus we find $c\geq{2}$ by Eq.~\eqref{eq: bound one species F} (a.k.a. the Fisher--KPP travelling wave speed).
For $m=1$, using the fact that $\Gamma(z+1)=z\Gamma(z)$, we see that $c\geq 1/\sqrt{2}$, and for the case $m = 2$, we find the explicit solution $c\geq \sqrt{(10+7\sqrt{7})/270}$. 
We highlight that this analysis also applies to the more general case where $D(u) = u^p$, and $f(u) = u^q(1-u)$, by defining $m = p + q$.
The next example illustrates a contrasting scenario: a pushed-wave regime where the variational principle becomes exact.

\subsection{A model with degenerate density-dependent diffusion and the Allee effect}

\begin{figure}[h!]
    \centering
    \includegraphics[width=0.8\linewidth]{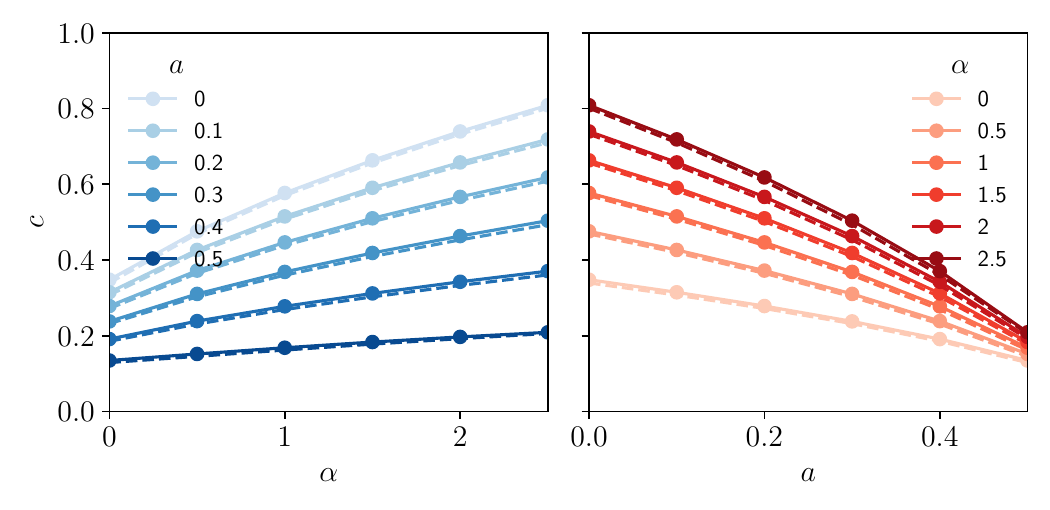}
    \caption{Plot of the numerically estimated minimum travelling wave speed (solid lines) and the minimum travelling wave speed of Eq.~\eqref{eq:allee} derived using the variational principle with $F(\beta)$ found in Eq.~\eqref{eq:f_allee} (dashed lines) as a function of $\alpha$ (left-hand side) and $a$ (right-hand side). The numerically estimated travelling wave speed is obtained by tracing the point $X(t)$ such that $u(X(t), t)=0.1$.}
    \label{fig:allee}
\end{figure}

In this example we consider a degenerate, parabolic reaction--diffusion equation, adapted from that studied in~\cite{sanchez1994travelling} to demonstrate how the variational approach can be used to estimate the travelling wave speed when the Allee effect is present~\cite{stephens1999allee}. 
The governing equation takes the form 
\begin{equation}
    \partial_t \rho = \partial_x ((\alpha \rho + \rho^2)\partial_x \rho) + \rho(1-\rho)(\rho-a), \label{eq:allee}
\end{equation}
where $D(\rho)=\alpha \rho + \rho^2$ and $f(\rho)=\rho(1-\rho)(\rho-a)$ for $\alpha>0$ and $a\in[0, 0.5]$. 
Employing the variational principle in Eq.~\eqref{eq: bound one species F}, we find
\begin{align}
    F(\beta)&=\dfrac{\beta}{120\Gamma(2-\beta)}[6(\alpha-a)\Gamma(4-\beta)+\Gamma(5-\beta)-30\alpha a\Gamma(3-\beta)] \nonumber \\
    &=\frac{\beta(2-\beta)}{120}\left[(3-\beta)(4-\beta+6(\alpha-a))-30\alpha a\right], \label{eq:f_allee}
\end{align}
where the recursion formula is used to simplify the Gamma functions.  
Employing Eq.~\eqref{eq:f_allee} in the expression for the minimum travelling wave speed given in Eq.~\eqref{eq: bound one species F}, we numerically solve for the supremum and find a very good prediction of the numerically observed travelling wave speed, which can be observed in Fig.~\ref{fig:allee}.

It is notable that the accuracy of the variational bound in this example is markedly higher than in the previous Porous--Fisher type case, and indeed becomes exact for certain parameter regimes.
This improvement is not accidental: it reflects key structural differences introduced by the Allee effect.
When $f(u)=\rho(1-\rho)(\rho-a)$, the reaction term no longer exhibits logistic behaviour and the classical marginal-stability mechanism (where the wave speed is determined solely by the linearisation near $u=0$) is replaced by a pushed-wave regime in which the entire interior of the wave profile contributes to speed selection. In contrast, for higher powers of the porous--Fisher equation, no admissible test function can be shown to attain equality in the variational principle.

In biological terms, the Allee effect suppresses growth at low densities, shifting the leading-edge dynamics away from the small $u$ asymptotics that the variational principle relies on, but also, in turn, reduces the accuracy of the variational bound.
As a consequence, the variational method instead here captures the dominant nonlinear contributions, producing a near-exact estimate across the full parameter ranges of $\alpha$ and $a$, as shown in Fig.~\ref{fig:allee}.

\subsection{Fisher--Stefan model}
In this example, the variational principle is adapted to study a well-known moving boundary problem, the Fisher--Stefan model, which is described by
\begin{align}
    \partial_t \rho & = \partial_x^2\,\rho+\rho(1-\rho)\,, \quad x < s(t) \label{eq:fisher_stefan}
    \\
    \rho & = 0,\quad \frac{\mathrm{d}s}{\mathrm{d}t} = -\kappa \partial_x \rho\,,\quad \mbox{on } x = s(t)\,, \label{eq:fisher_stefan_bc}
\end{align}
where $\kappa\geq 0$ is a parameter that controls the speed of the moving boundary.

By fixing the moving boundary at $z=0$, we can use the travelling wave ansatz, $z = x-s(t) = x-ct$ with $\rho(x,t)=u(z)$ to obtain
\begin{align}
    u'' & + cu' + u(1-u)  = 0 \label{eq:fs-tw}
    \\
   u(0) & = 0,\quad u'(0) = -\frac{c}{\kappa}\,. \nonumber
\end{align}
As before, expressing $v(u) = -\mathrm{d}u/\mathrm{d}z$ as a function of $u$, and integrating Eq.~\eqref{eq:fs-tw} against a test function, $\varphi(u)$, we obtain
\begin{equation*}
    \int_0^1 \left(-\frac{1}{2}\phi(u) v^2(u)+cv(u)\varphi(u)\right)\mathrm{d}u  + \dfrac{c^2 \varphi(0)}{2\kappa^2}= \int_0^1 u(1-u)\varphi(u)\,\mathrm{d}u\,,
\end{equation*}
where we used that $v(0) = c/\kappa$. 
Using the same argument as previously shown on the second order polynomial in $v$, we obtain the bound
\begin{equation*}
   \frac{1}{2} c^2 \geq \frac{\int_0^1 u(1-u)\varphi(u)\,\mathrm{d}u}{\varphi(0)/\kappa^{2}+\int_0^1 \varphi^2(u)/\phi(u)\,\mathrm{d}u}\,,
\end{equation*}
with equality if
\begin{equation*}
    v = \frac{c\varphi(u)}{\phi(u)}\,\Longleftrightarrow \varphi(u) \propto \exp\left\{-\int_0^u\frac{c\,\mathrm{d}q}{v(q)}\right\} .
\end{equation*}
Evaluating at $u = 0$, we also obtain the compatibility condition
\begin{equation*}
    \varphi(0) = \phi(0)/\kappa\,,
\end{equation*}
which means that $\varphi(u)= e^{-\kappa u} \Tilde{\varphi}(u)$ with $ \Tilde{\varphi}'(0) = 0$. 
Without loss of generality we assume $\varphi(0) = \Tilde{\varphi}(0) = 1$. 

\begin{figure}[htbp]
    \centering
    \includegraphics[width=.95\linewidth]{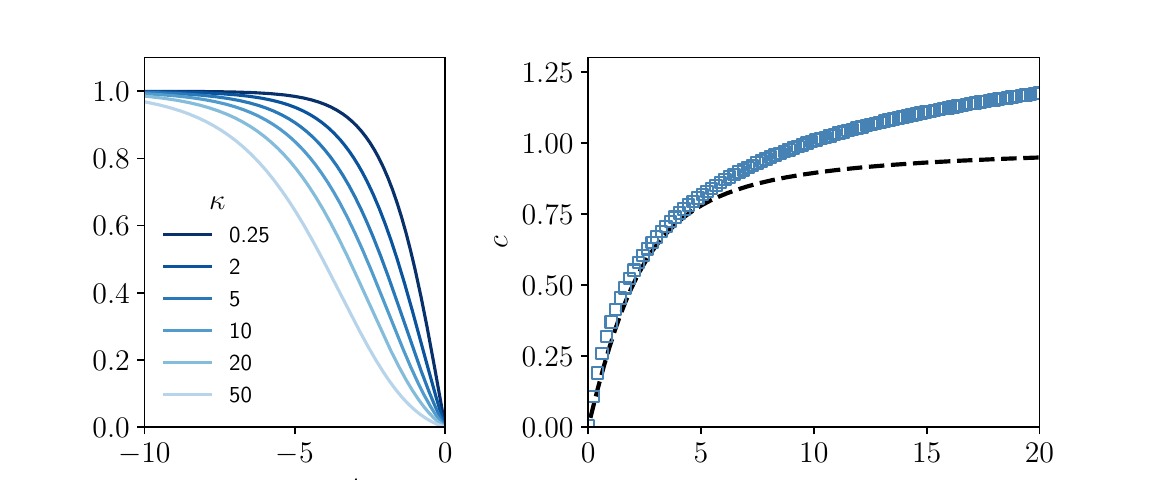}
    \caption{(left) Shape of travelling wave profiles for the Fisher-Stefan model given by Eqs.~\eqref{eq:fisher_stefan} and \eqref{eq:fisher_stefan_bc} depending on the parameter $\kappa$. (right) Plot of the numerically simulated minimum travelling wave speed (squares) and analytically calculated travelling wave speed using the variational principle in Eq.~\eqref{eq: bound fisher stefan simple} (black dashed lines) as a function of $\kappa$. The numerically estimated travelling wave speed is obtained by tracing the point $X(t)$ such that $u(X(t), t)=0.1$.}
    \label{fig:fisher-stefan}
\end{figure}

The simplest test function we could employ that satisfies these bounds and the compatibility condition is $\varphi(u) = e^{-\kappa u}$.
Using this, we find
\begin{equation}
    \dfrac{1}{2}c^2\geq \dfrac{\int_0^1 u(1-u) e^{-\kappa u}\,\mathrm{d}u}{1/\kappa^2 + \int_0^1e^{-\kappa u}/\kappa\,\mathrm{d}u} = \frac{\kappa-2+e^{-\kappa}(2+\kappa)}{\kappa(2+e^{-\kappa})}\,,\,\label{eq: bound fisher stefan simple}
\end{equation}
which gives the correct behaviour~\cite{el2021invading} as $\kappa \rightarrow 0^+$: 
\begin{equation*}
    c\sim \frac{\kappa}{\sqrt{3}}\,.
\end{equation*}
Eq.~\eqref{eq: bound fisher stefan simple} therefore provides a lower bound for the travelling wave speed that is valid for all $\kappa > 0$; however, the quality of this bound degenerates as $\kappa\rightarrow +\infty$ (see Fig.~\ref{fig:fisher-stefan}, where it is notable that the bound tends to one in the limit $\kappa\to\infty$ while the speed should actually tend to two~\cite{mccue2021exact}).


The degeneration of the bound in the limit $\kappa\to\infty$ is not just a shortcoming of the chosen test function, but reflects a fundamental limitation of the variational formulation in this regime. 
As $\kappa$ increases, the influence of the Stefan boundary condition weakens and the Fisher--Stefan model approaches the classical Fisher--KPP equation, for which the travelling wave speed is selected by linear marginal stability at the leading edge. The variational principle employed here encodes the moving boundary condition through constraints on admissible test functions at $u=0$. 
In the limit $\kappa\to\infty$, these constraints degenerate, and the class of admissible test functions no longer captures the linearised dynamics governing speed selection. As a consequence, the resulting bound saturates at a value strictly below the true wave speed $c=2$.

This observation highlights that, while the variational approach yields a rigorous and informative lower bound for all $\kappa>0$, it is most effective when wave speed selection is governed by nonlinear or boundary effects. 
Recovering the Fisher--KPP speed in the large $\kappa$ limit using these techniques would require a reformulation of the variational principle that explicitly incorporates the linear leading edge behaviour in this asymptotic regime, which lies beyond the scope of the present work.

\section{Two-species models: a leading order variational approximation}\label{sec:MS}
Although the full optimal control problem derived in Sec.~\ref{sec:OC} is exact, solving the associated optimality system analytically is generally intractable for models involving two or more species. 
Instead, in this section, we demonstrate how we can exploit the optimal control approach to find a leading order variational approximation that captures the dominant wave speed selection mechanism in a weak coupling regime between the two species.
The lower bound for the minimum travelling wave speed found using this variational approach improves the bound found using linear stability analysis, and thus is  particularly useful in specific parameter regimes where further analytical progress is made possible. 

In particular, we restrict our attention to a particularly relevant class of models where $\rho_2$ denotes a density field which gets degraded by the first species, such that  $g(\rho_1,\rho_2) = -\kappa\rho_1\rho_2$ with $\rho_2(+\infty) = \nu$. 
Solving for $u_2$, we find 
\begin{equation}
   u_2(u_1) = \nu\exp\left\{-\frac{\kappa}{c}\int_0^{u_1}\frac{q\,\mathrm{d}q}{v(q)}\right\}.\label{eq: u_2 solution}
\end{equation}
In this case, weak coupling corresponds to small values of $\kappa$ and $\nu$, such that the secondary species remains uniformly small across the wavefront. 
Indeed, Eq.~\eqref{eq: u_2 solution} shows that $u_2(u_1)=\mathcal{O}(\nu)$ uniformly on $[0,1]$.
We therefore introduce the small parameter
$\epsilon = {\kappa\nu}/{c},$
and seek an approximate solution valid as $\epsilon\to0$.

We start by noting that for two-species systems, the optimal control $v^*$ can be found by solving the equations
\begin{align}
\label{eq: oc 2species v}   0 &  = -\phi D^2 v^3+cD\varphi v^2 + \frac{\kappa}{c}{\omega u_1 u_2} \,,
   \\
   \frac{\mathrm{d}u_2}{\mathrm{d}u_1} & = -\frac{\kappa}{c}\frac{u_1 u_2}{v}\,, \label{eq:u2_u1}
   \\
   \frac{\mathrm{d}\omega}{\mathrm{d}u_1} & = \left[\phi D v^2-cv\varphi\right]\partial_2 D + \partial_2 (Df)\varphi+\frac{\kappa}{c}\frac{u_1 \omega}{v}\,, \label{eq:omega_u1}
\end{align}
with $u_2(0) = \nu$, $u_2(1) = 0$, and $\partial_2 = \partial/\partial_{u_2}$, $D = D(u_1,u_2)$ and $f=f(u_1,u_2)$, for notational simplicity. 

Motivated by the variational principle for single--species models, we seek a solution of the form 
$    v= {c\varphi}/{\phi}D+\delta v$
where $\delta v=\mathcal{O}(\epsilon)$ uniformly in $[0,1].$
First, we observe that the adjoint enters the stationarity condition (Eq.~\eqref{eq: oc 2species v}) only through the product $u_2\omega$. 
Differentiating this product and using Eqs.~\eqref{eq:u2_u1} and~\eqref{eq:omega_u1} yields 
\begin{equation*}
    \frac{\mathrm{d}}{\mathrm{d}u_1}(u_2\omega) = u_2\left[\phi D v^2-cv\varphi\right]\partial_2 D + u_2\partial_2 (Df)\varphi\,,
\end{equation*}
and hence
\begin{equation*}
    u_2(u_1)\omega(u_1) = -\int_{u_1}^1 \left(u_2\left[\varphi D v^2-cv\varphi\right]\partial_2 D + u_2\partial_2 (Df)\varphi\right)\,\mathrm{d}q\,.
\end{equation*}

Since the secondary species satisfies $0\le u_2(u_1)\le \nu$ uniformly on $[0,1]$, and admissible test functions obey $\varphi(u_1)\lesssim u_1^{-\beta}$ with $\beta\in(0,2)$, the integral representation for the adjoint contribution $u_2\omega$ remains uniformly bounded on $[0,1]$.
As a consequence, $u_2\omega=\mathcal{O}(\nu)$ uniformly, and the adjoint term in the stationarity condition (Eq.~\eqref{eq: oc 2species v}) enters only at order $\mathcal{O}(\kappa\nu/c)=\mathcal{O}(\epsilon)$.

Under mild regularity assumptions on the reaction and diffusion functions, and standard admissibility conditions on the test function, we show in Appendix~\ref{app:reg} that the ansatz $v = c\varphi/\phi + \delta v$ satisfies the full optimality condition (Eq.~\eqref{eq: oc 2species v}) up to order $\mathcal{O}(\epsilon)$.
As such, at leading order in the weak coupling parameter $\epsilon$, the optimal control formulation coincides with that of the one-species variational principle, given by
\begin{equation}
    \frac{1}{2}c^2 \geq \frac{\int_0^1 D(u_1, u_2(u_1))f(u_1, u_2(u_1))\varphi(u_1)\,\mathrm{d}u_1}{\int_0^1 \varphi(u_1)^2/\phi(u_1)\,\mathrm{d}u_1},\quad\mbox{with } \phi = - \varphi'\,, \label{eq:VP_MS}
\end{equation}
where $u_2(u_1)$ is given by Eq.~\eqref{eq: u_2 solution} and equality is only achieved when 
\begin{equation}
    v^* = \frac{c\varphi(u_1)}{D(u_1,u_2)\phi(u_1)}\,.
    \label{eq: ode constraint}
\end{equation}

This leading order reduction provides a practical variational bound for multispecies systems. 
In the following subsections, we apply this framework to biologically motivated models and demonstrate that it captures nonlinear wave speed selection beyond linear theory.

\subsection{Models of cell invasion into extracellular matrix}
We now consider a class of models describing collective cell migration into the extracellular matrix (ECM) which take the following form 
\begin{align}
    \partial_t \rho_1 & = \partial_x\left((1-\rho_2)\partial_x \rho_1\right) + f(\rho_1, \rho_2)\,,
    \label{eq:ECM_u}
    \\
    \partial_t \rho_2 & =   -\kappa\rho_1\rho_2\,,
    \label{eq: ECM_m}
\end{align}
where $f(0, \rho_2)=0$ and $f(1, \rho_2)=0.$
The boundary conditions are given by
\[
\begin{cases}
\rho_1(x,t) \to 1,\quad \rho_2(x,t)\to 0, &\text{as } x \to -\infty,  \\
\rho_1(x,t) \to 0,\quad \rho_2(x,t)\to \nu, &\text{as } x \to +\infty,
\end{cases}
\]
for $0\leq\nu\leq1.$

We consider two different reaction terms.
Firstly, $f_C(\rho_1,\rho_2) = f_C(\rho_1) = \rho_1(1-\rho_1)$ as described in~\cite{colson2021travelling}.
Secondly, we consider $f_B(\rho_1,\rho_2) = \rho_1(1-\rho_1-\rho_2)$ as detailed in~\cite{browning2019bayesian, el2021travelling}. 
For both models, numerical simulations show that the system exhibits travelling waves that attain the linearly predicted wave speed ($2\sqrt{1-\nu}$ when $f=f_C$, and $2(1-\nu)$ when $f=f_B$) for small values of the ECM degradation rate, $\kappa \leq \kappa^*(\nu)$, only, where $\kappa^*(\nu)$ is the critical ECM degradation rate such that the linearly predicted and numerically observed travelling wave speeds coincide for $\kappa \leq \kappa^*(\nu)$.
For sufficiently large ECM degradation rates $\kappa > \kappa^*(\nu)$, this model exhibits nonlinear wave speed selection.  
The associated numerically estimated wave speed show significant deviations from the linearly selected travelling wave speed, which underestimates the invasion speed.
Furthermore, the critical value of the degradation rate, $\kappa^*(\nu)$, seems to follow a monotonically increasing relationship with $\nu$ as $\nu\rightarrow 1$~\cite{crossley2023travelling}.

We now want to apply the leading order variational approximation from Sec.~\ref{sec:MS}.
Employing the same family of test functions as in Sec.~\ref{1-vp}, ($\varphi_\beta (u_1) = ((1-u_1)/u_1)^\beta$ for $\beta \in[0,2)$), we have 
\begin{align}
    \dfrac{1}{2}c^2\geq \sup_{\beta\in[0,2)} G(\beta),\qquad G(\beta) = \beta\,\frac{M(\beta)}{\mathrm{B}(2-\beta, 2+\beta)}\,, \label{eq:N_vp}
\end{align}
where
\begin{equation}
    M(\beta) = \int_0^1 D(u_1, u_2(u_1))f(u_1, u_2(u_1))\left(\dfrac{1-u_1}{u_1}\right)^\beta\,\mathrm{d}u_1. \label{eq:N}
\end{equation}
For $D(u_1, u_2(u_1))=(1-u_2(u_1))$ and $f_C(u_1) = u_1(1-u_1)$, we find
\begin{align}
    M_C(\beta)&=\int_0^1 \dfrac{u^{1-\beta}(1-u)^{1+\beta}}{\frac{\nu}{1-\nu}(1-u)^{{\kappa\beta}/{c^2}+1}}\mathrm{d}u. \label{eq:N1}
\end{align}
However, for $f_B(u_1,u_2) = u_1(1-u_1-u_2)$ with $\gamma=\nu/(1-\nu)$ we have 
\begin{align}
    M_B(\beta) &= \int_0^1 (1-u_2(u_1))(1-u_1-u_2(u_1))\,u_1^{1-\beta}(1-u_1)^{\beta}\,\mathrm{d}u_1 \nonumber \\
    &= \int_0^1 \dfrac{u_1^{1-\beta}(1-u_1)^{\beta}-u_1^{2-\beta}(1-u_1)^\beta-\gamma(1-u_1)^{1-\beta}(1-u_1)^{\beta+\eta}}{(1+\gamma(1-u_1)^\eta)^2}. \label{eq:N2}
\end{align}

\begin{figure}[htbp]
    \centering
\includegraphics[width=\linewidth]{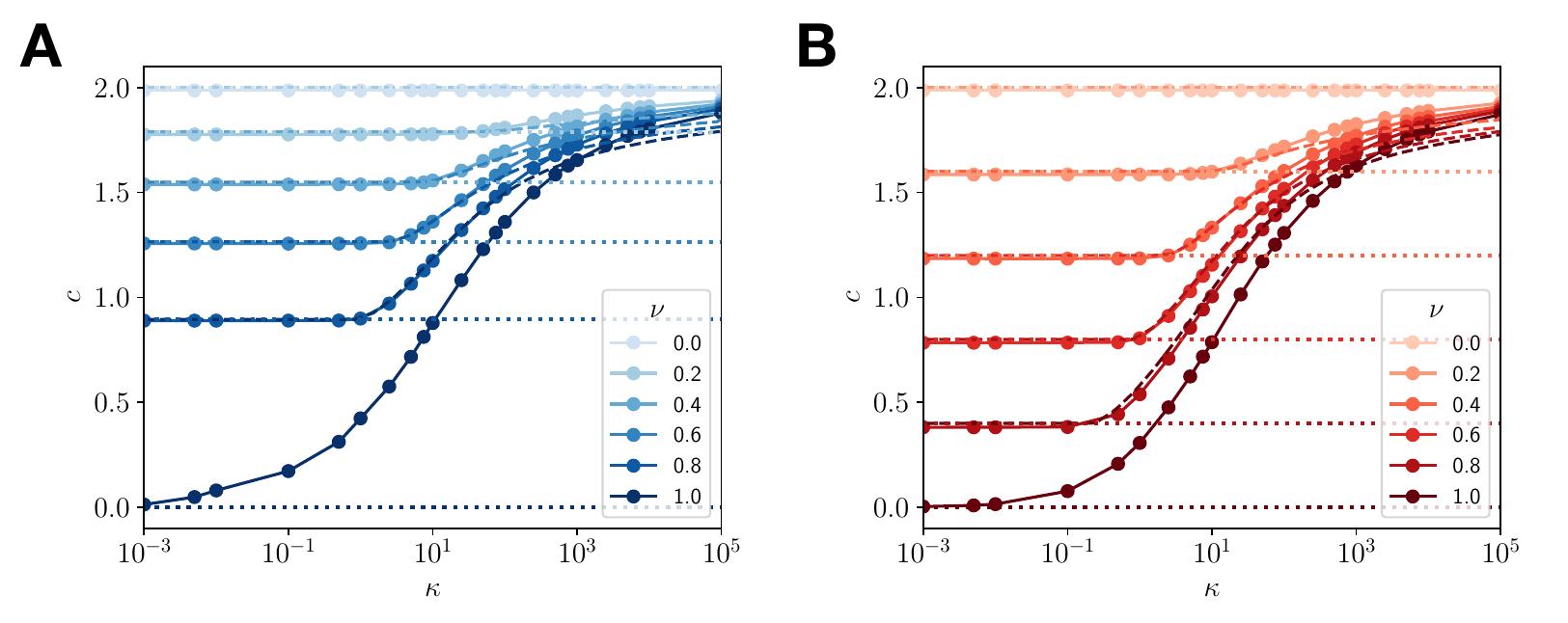}
    \caption{
    Plot of numerically simulated minimum travelling wave speed (solid lines) and analytically calculated travelling wave speed of Eqs.~\eqref{eq:ECM_u} and \eqref{eq: ECM_m} using linear theory (dotted lines) and the variational principle (dashed lines) with $M_C(\beta)$ as calculated in Eq.~\eqref{eq:N1} (left-hand side) and $M_B(\beta)$ as calculated in Eq.~\eqref{eq:N2} (right-hand side) for various initial ECM densities $\nu\geq0$. 
    The numerically estimated travelling wave speed is obtained by tracing the point $X(t)$ such that $u_1(X(t), t)=0.1$.}
    \label{fig:colson}
\end{figure}

Substituting Eqs.~\eqref{eq:N1} and~\eqref{eq:N2} into Eq.~\eqref{eq:VP_MS} provides a lower bound on the travelling wave speed that, in each case, can be solved for numerically.
Fig.~\ref{fig:colson}A shows the result of numerically solving Eq.~\eqref{eq:VP_MS} subject to Eq.~\eqref{eq:N1} (a.k.a. when $f(u_1, u_2)=f_C(u_1)$). Similarly,  Fig.~\ref{fig:colson}B shows the result of numerically solving Eq.~\eqref{eq:VP_MS} subject to Eq.~\eqref{eq:N2} (a.k.a. when $f(u_1, u_2)=f_B(u_1, u_2)$) for a lower bound on the travelling wave speed, $c$, which shows a vast improvement when compared to the linearly predicted travelling wave speeds.
To the best of our knowledge, the linear travelling wave speed is the current best estimate of the speed of invasion in the literature for both of the models presented in this section, so the prediction derived in this work using the variational approach drastically improves the lower bound on the wave speed in both cases.

\subsection{A model of cell invasion with cell differentiation} 
\begin{figure}
    \centering
    \includegraphics[width=0.98\linewidth]{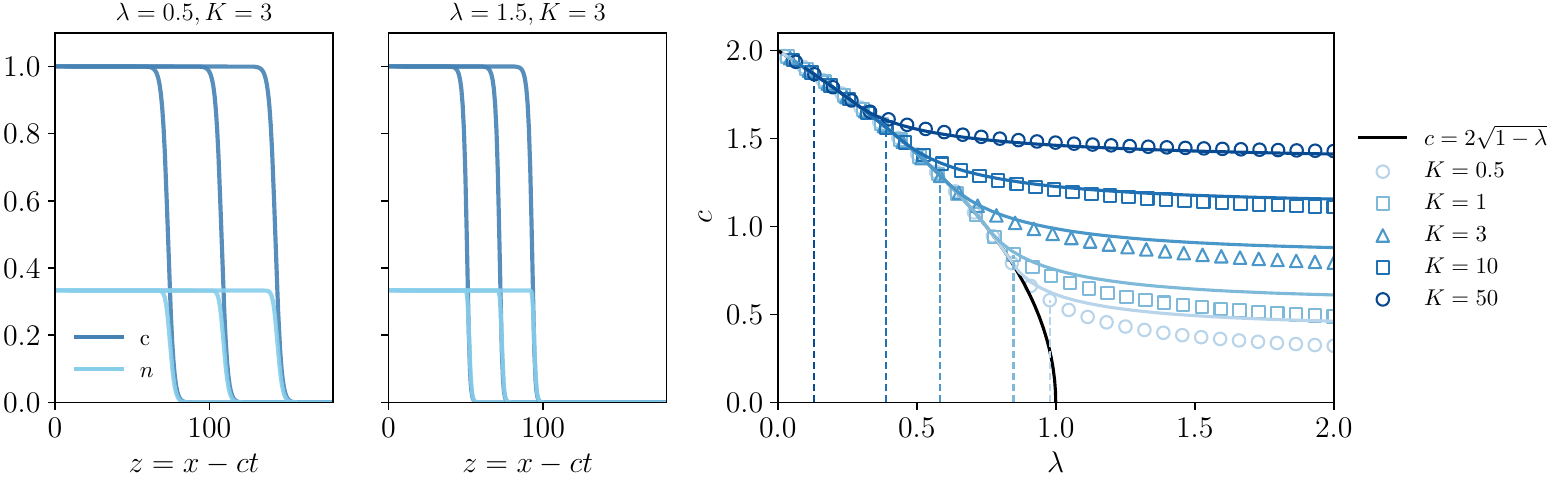}
    \caption{
(left) Numerical simulations of the time-dependent cell differentiation model (Eqs.~\eqref{eq: landman a}  and~\eqref{eq: landman b}) exhibiting travelling wave solutions. Solutions plotted at $t = 50,\,75,\,100$ for different parameters. (right) Wave speed $c$ as a function of $\lambda$ for various values of $K$.
Solid lines correspond to predictions obtained from the approximate variational principle, while square markers indicate numerical solutions of the full PDE system.
Dashed vertical lines mark the threshold $\lambda \kappa / c = 1$, beyond which the weak coupling approximation is no longer valid. The black curve shows the linear theory prediction $c = 2\sqrt{1-\lambda}$.
}
    \label{fig: fig5}
\end{figure}

We now consider a two--species model of cell invasion first introduced in~\cite{trewenack2009traveling}, which has recently attracted renewed interest due to its transition from linear to nonlinear wave speed selection~\cite{yue2025traveling}. The model describes a population of invasive cells, $\mathrm{c}(x,t)$, that migrate via diffusion and proliferate logistically, while differentiating at rate $\lambda>0$ into a population of non--motile, non--proliferative cells, $n(x,t)$---see Fig.~\ref{fig: fig5} for numerical solutions. The governing equations are
\begin{align}
    \partial_t \mathrm{c} & = \partial_x^2 \mathrm{c} + \mathrm{c}(1-\mathrm{c}) - \lambda \mathrm{c}(1-K n)\,,\label{eq: landman a}
    \\
    \partial_t n & = \lambda \mathrm{c}(1-K n)\,,\label{eq: landman b}
\end{align}
for $K > 0$ modulating the differentiation rate. The steady states of this system are $(\mathrm{c},n) = (0,n^*)$ and $(1,1/K)$, for any $0\leq n^*\leq 1/K$, where we choose $n^* = 0$ as this corresponds to the more natural case of having no differentiated cells initially~\cite{trewenack2009traveling}. 
This model can be reduced to a form similar to that in Eqs.~\eqref{eq:ECM_u} and~\eqref{eq: ECM_m} by setting $\rho_2 = 1-K n$, $\rho_1 = \mathrm{c}$. 
This gives
\begin{align}
    \partial_t \rho_1 & = \partial_x^2\, \rho_1 + \rho_1\left(1-\rho_1-\lambda \rho_2\right)\,, \label{eq: landman u ctrl}
    \\
    \partial_t \rho_2 & = - \kappa \rho_1\rho_2\,, \label{eq: landman m ctrl}
\end{align}
where $\kappa = \lambda K$. We assume that $\lambda K \ll 1$ so that we can use the leading order variational principle we derived in Sec.~\ref{sec:MS}. 
The boundary conditions of the system are given by
\[
\begin{cases}
\rho_1(x,t) \to 1\quad \rho_2(x,t)\to 0, &\text{as } x \to -\infty\,,  \\
\rho_1(x,t) \to 0\quad \rho_2(x,t)\to 1, &\text{as } x \to +\infty\,.
\end{cases}
\]
Setting $\rho_1(x,t) = u_1(x-ct)$ and $\rho_2(x,t) = u_2(x-ct)$, the linear marginal stability condition in this case gives $c \geq c_L = 2\sqrt{1-\lambda}$ which only applies for $\lambda \leq 1$.
Employing the leading order variational principle using Eq.~\eqref{eq:N}, we calculate
    \begin{align*}
        M(\beta) & = \int_0^1 u_1^{1-\beta}(1-u_1)^\beta\left[(1-u_1)-\lambda(1-u_1)^{\beta\kappa/ c^2}\right]\,\mathrm{d}u_1
        \\
        & =  \mathrm{B}(2-\beta, 2+ \beta) -\lambda B\left(2-\beta, 1+\beta+\frac{\beta\kappa}{c^2}\right)\,.
    \end{align*}
Hence, Eq.~\eqref{eq:N_vp} gives
    \begin{align}
        G(\beta) = \beta\left[1-6\lambda\frac{\Gamma\left(1+\beta+\frac{\beta\kappa}{ c^2}\right)}{\Gamma\left(3+\frac{\beta\kappa}{ c^2}\right)\Gamma(2+\beta)}\right]\,.
        \label{eq: Gc landman model}
    \end{align}
We note that $G(2) = 2(1-\lambda)$, giving $c^2 \geq 4(1-\lambda)$. This is in agreement with the travelling wave speed predicted by linear theory. In particular, our approach recovers, and in fact improves upon, the marginal bound obtained from linear stability analysis.

Next, we use Eq.~\eqref{eq: Gc landman model} to obtain the sharpest bound on the wave speed in the form ${c^2}/{2} \geq \sup_{\beta\in(0,2)} G(\beta)$. In Fig.~\ref{fig: fig5} we also indicate the parameter regimes in which the weak coupling approximation breaks down, namely when ${\lambda \kappa}/{c} \sim 1$. Overall, we observe that the approximate variational principle provides an accurate prediction of the wave speed. Agreement is excellent within the weak coupling regime, and although deviations become visible as this regime is exited, the resulting error remains comparatively small.

It is particularly interesting that the approximation remains accurate even for large values of $K$. This motivates an investigation of the cell differentiation model in the regime where $\kappa = \lambda K$ is large but $\lambda$ remains small. In this case the adjoint variable $\omega$ satisfies
\begin{equation*}
    \omega(u_1)u_2(u_1) = \lambda\int_{u_1}^1q\varphi(q)u_2(q)\,\mathrm{d}q\,.
\end{equation*}
It is easy to show that when $\kappa\rightarrow +\infty$ and $\lambda\ll 1$ the ansatz $v^* =c\varphi/\phi$ is the leading order solution to the optimal control formulation. Indeed, using this choice of $v$ to solve for $u_2$, with $\varphi(u_1) = ((1-u_1)/u_1)^\beta$, and expanding for large $\kappa$, we obtain
\begin{align*}
     \omega(u_1)u_2(u_1) = \lambda\int_{u_1}^1 q^{1-\beta}(1-q)^{\beta + \kappa\beta/c^2}\,\mathrm{d}q \sim u_1^{1-\beta}\frac{(1-u_1)^{\kappa\beta/c^2}}{\beta + \kappa\beta/c^2}\,.
\end{align*}
In particular, as $\kappa\rightarrow +\infty$, we obtain
\begin{equation}
    \frac{\kappa}{c}\omega u_1 u_2 = \mathcal{O}(c\lambda)\,,
\end{equation}
so the weak coupling assumption is in fact satisfied if $\lambda c$ is small, which explains the good agreement observed in Fig.~\ref{fig: fig5}.

\subsection{Extension to models with taxis terms}
Benguria et al. have, more recently, considered the addition of advection terms in the governing equations and their impact on the resulting travelling wave speed~\cite{BENGURIA2022112668}.
As is natural, we therefore present a simple extension of the work presented here that is applicable to populations undergoing taxis, or following a gradient of external signals, for example, chemicals or substrates---\textit{e.g.}, chemotaxis, haptotaxis, electrotaxis.
The models of interest will take the form 
\begin{align}
    \partial_t \rho_1 & = \partial_x\left(D(\rho_1,\rho_2)\partial_x \rho_1-\chi(\rho_1, \rho_2)\partial_x\rho_2\right) + f(\rho_1, \rho_2)\,,
    \label{eq: reaction diffusion two species a_chem}
    \\
    \partial_t \rho_2 & =   g(\rho_1, \rho_2)\,,
    \label{eq: reaction diffusion two species b_chem}
\end{align}
where the boundary conditions of the system are given by
\[
\begin{cases}
\rho_1(x,t) \to 0\quad \rho_2(x,t)\to \mu, &\text{as } x \to -\infty,  \\
\rho_1(x,t) \to 1\quad \rho_2(x,t)\to \nu, &\text{as } x \to +\infty.
\end{cases}
\]

If we have that $D(\rho_1,\rho_2)\geq 0$ is the density-dependent diffusivity, $\chi(\rho_1,\rho_2)\geq 0$ is the tactic sensitivity, and $f(\rho_1,\rho_2), g(\rho_1,\rho_2)$ are reaction terms,  we can then look for a monotonically decreasing function of $z = x-ct$, defining $\rho_1(x,t) = u_1(x-ct)$ and $\rho_2(x,t) = u_2(x-ct)$, such that 
\[
\begin{cases}
u_1(z) \to 0 \quad u_2(z) \to \mu \quad \text{as } z \to -\infty, \\
u_1(z) \to1 \quad u_2(z) \to \nu \quad \text{as } z \to \infty.
\end{cases}
\]
If we consider a weak coupling with the external cue, we can again follow the variational principal approach detailed earlier in Sec.~\ref{sec:MS}. Namely, setting $\mathrm{d}u_1/\mathrm{d}z  = -v$, we have
\begin{align*}
    \dfrac{\mathrm{d}u_1}{\mathrm{d}z} &= -v, \\
    \dfrac{\mathrm{d}}{\mathrm{d}z}\bigg(D(u_1,u_2)v+\chi(u_1, u_2)\frac{\mathrm{d}u_2}{\mathrm{d}z}\bigg) &= f(u_1,u_2)-cv, \\
    \dfrac{\mathrm{d}u_2}{\mathrm{d}z}  &=-\dfrac{1}{c}g(u_1,u_2),
\end{align*}
Replacing all instances of $D(u_1,u_2)v$ by $D(u_1,u_2)v+\chi(u_1, u_2){\mathrm{d}u_2}/{\mathrm{d}z}$ in the calculations, we can follow our previous working to show an approximate variational principle for chemotactic style models. 
Ultimately, we find
\begin{equation*}
    \frac{1}{2}c^2 \geq \frac{\int_0^1 f(u_1, u_2(u_1))\varphi(u_1)\bigg(D(u_1,u_2)-\chi(u_1, u_2)\frac{\mathrm{d}u_2}{\mathrm{d}u_1}\bigg)\,\mathrm{d}u_1}{\int_0^1 \varphi(u_1)^2/\phi(u_1)\,\mathrm{d}u_1},\quad\mbox{with } \phi = - \varphi'\,,
\end{equation*}
as long as 
\begin{equation}
    \frac{\mathrm{d}u_2}{\mathrm{d}u_1} = \frac{\phi(u_1)\bigg(D(u_1,u_2)-\chi(u_1, u_2)\frac{\mathrm{d}u_2}{\mathrm{d}u_1}\bigg)g(u_1,u_2)}{c^2\varphi(u_1)}.
    \label{eq: ode constraint_chem}
\end{equation}
Due to lack of examples in the literature, this method is not applied to any models in this work but instead detailed to demonstrate the utility and wider applicability of the approach to these models, where they exist.

\section{Discussion}\label{sec:conc}

In this article, we have developed an optimal control framework for characterising the minimal travelling wave speed in systems of reaction--diffusion equations with nonlinear diffusion or reaction terms admitting monotonic fronts. 
This approach yields a variational characterisation of the minimal speed that recovers classical results for single--species equations~\cite{benguria1994validity, benguria1996speed, BENGURIA2022112668, stokes2024speed} and, in a weak coupling regime, extends them naturally to coupled systems. 
Within this framework, we also derive a criterion for the onset of nonlinear wave speed selection, and apply the associated variational principle to several biologically motivated models to elucidate parameter dependence.

From an optimal control perspective, the classical single--species variational principle emerges as a special, explicitly solvable case of a more general control problem: the travelling wave profile serves as the state variable, the local flux as the control, and the wave speed arises as the optimal value satisfying a stationarity condition. 
This interpretation clarifies the structure underlying existing variational arguments and provides a systematic route to generalising wave speed selection theory to multi-species systems, where closed-form variational formulations are typically unavailable.

Despite the benefits of the variational characterisation, further research is needed to address the limitations of the methods. 
The analysis assumes the existence of monotonic travelling waves with constant speed, enabling the reduction of the PDEs to a system of ordinary differential equations in the travelling wave coordinate. 
Moreover, the multi-species variational approximation relies on sufficiently weak coupling so that the adjoint contributions in the optimality system may be treated as higher-order corrections only. 
In regimes with strong coupling, sharp internal layers, or strong dependence of the diffusivity of the motile species on the additional species, the resulting lower bounds are not expected to be sharp and may only capture qualitative parameter dependence. 
Nonetheless, the framework presented here provides improved estimates on the lower bound for the travelling wave speed together with a unified and extensible foundation for studying nonlinear wave speed selection in general in systems of reaction--diffusion equations.

\appendix
\section{Uniform boundedness of the adjoint contribution}\label{app:reg}
In order to show that if $v = c\varphi/D\phi + \delta v$, where $\delta v = \mathcal{O}(\epsilon)$, then $v$ satisfies Eq.~\eqref{eq: oc 2species v} uniformly in $u_1\in[0,1]$ up to order $\mathcal{O}(\epsilon)$, we make the following mild assumptions on the reaction and diffusivity of the first species, and on the test function. 

We assume that there exist constants $a_1, a_2, a_3, a_4> 0$ such that
\begin{equation}
\max_{0\leq u_1\leq 1,\,0\leq u_2\leq \nu}D, |\partial_2 D|, |\partial_2(Df)| \leq a_1,\quad D \geq a_2 > 0;
\end{equation}
and
\begin{equation*}
   a_3(1-u_1)^{\beta} \leq\varphi\leq a_4 u^{-\beta},\qquad a_4(1-u_1)^{\beta-1}\leq\phi\leq a_4 u^{-\beta-1}\,,
\end{equation*}
for $\beta\in(0,2)$. 

Consider the final term, involving $\omega$, in Eq.~\eqref{eq: oc 2species v}. 
Using the boundedness assumptions, we find that
\begin{align*}
    \left|u_1 u_2 \omega\right|\leq  C\left(\nu c^2 u_1\int_{u_1}^1 q^2(1-q)^2\,\mathrm{d}q + \nu u_1\int_{u_1}^1 u^{-\beta}\right) +\mathcal{O}(\epsilon)\,.
\end{align*}
And from here we obtain 
\begin{equation*}
\sup_{\beta\in(0,2)}\max_{u_1\in[0,1]}|u_1 u_2 \omega| \leq \nu C + \mathcal{O}(\epsilon)\,.
\end{equation*}
where $C =C(a_1, a_2, a_3, a_4, c) > 0$. Thus $v^* = c\varphi/D\phi$ solves the optimal control problem, uniformly in $u_1\in[0,1]$, up to order $\mathcal{O}(\epsilon)$.

\section*{Acknowledgments}
The authors would like to thank Mat Simpson and Philip Maini for helpful discussions.  For the purpose of Open Access, the authors have applied a CC BY public copyright licence to any Author Accepted Manuscript (AAM) version arising from this submission. 
R.M.C. would like to thank the Engineering and Physical Sciences Research Council (EP/T517811/1) for funding.
C.F. acknowledges support from a Hooke Research Fellowship.
R.E.B. acknowledges support of a grant from the Simons Foundation (MP-SIP-00001828).

\end{document}